\documentclass[11pt]{amsart}%
\usepackage{amssymb,amsmath,amsfonts,latexsym,amsthm,geometry,graphicx}
\usepackage{amsmath}
\usepackage{amsfonts}
\usepackage{amssymb}
\usepackage{graphicx}%
\providecommand{\U}[1]{\protect\rule{.1in}{.1in}}
\providecommand{\U}[1]{\protect\rule{.1in}{.1in}}
\providecommand{\U}[1]{\protect\rule{.1in}{.1in}}
\newtheorem{theorem}{Theorem}[section]

\theoremstyle{definition}

\newtheorem{remark}[theorem]{Remark}

\begin{document}
\title[On the constants of the multilinear Hardy--Littlewood
inequality ]{A short communication on the constants of the multilinear Hardy--Littlewood inequality}
\author[D. Pellegrino ]{Daniel Pellegrino}
\address{Departamento de Matem\'{a}tica \\
\indent Universidade Federal da Para\'{\i}ba \\
\indent 58.051-900 - Jo\~{a}o Pessoa, Brazil.}
\email{pellegrino@pq.cnpq.br and dmpellegrino@gmail.com}
\thanks{The author is supported by CNPq.}
\subjclass[2010]{11Y60, 47H60.}
\keywords{Bohnenblust--Hille inequality, Hardy--Littlewood inequality}
\maketitle

\begin{abstract}
It was recently proved that for $p>2m^{3}-4m^{2}+2m$ the constants of the
Hardy--Littlewood inequality for $m$-linear forms on $\ell_{p}$-spaces are
less than or equal to the best known estimates of respective constants of the
Bohnenblust--Hille inequality. In this note we obtain upper bounds for
opposite side, i.e., the constants when $2m\leq p\leq2m^{3}-4m^{2}+2m.$ For
these values of $p$ our result improves previous estimates from 2014 of Araujo
\textit{et al}. for all $m\geq3.$

\end{abstract}

\section{Introduction}

Let $\mathbb{K}$ be the real scalar field $\mathbb{R}$ or the complex scalar
field $\mathbb{C}$. For any integer $m\geq2$, the Hardy--Littlewood inequality
(see \cite{alb, dimant, hardy,pra}) for $m$-linear forms asserts that for
$2m\leq p\leq\infty$ there exists a constant $C_{m,p}^{\mathbb{K}}\geq1$ such
that, for all $m$--linear forms $T:\ell_{p}^{n}\times\cdots\times\ell_{p}%
^{n}\rightarrow\mathbb{K}$ and all positive integers $n$,
\begin{equation}
\left(  \sum_{j_{1},...,j_{m}=1}^{n}\left\vert T(e_{j_{1}},...,e_{j_{m}%
})\right\vert ^{\frac{2mp}{mp+p-2m}}\right)  ^{\frac{mp+p-2m}{2mp}}\leq
C_{m,p}^{\mathbb{K}}\left\Vert T\right\Vert . \label{i99}%
\end{equation}
Moreover, the exponents $\frac{2mp}{mp+p-2m}$ are optimal. The case $p=\infty$
recovers the famous Bohnenblust--Hille inequality (see \cite{bh}):

\textbf{Theorem }(Bohnenblust--Hille inequality). There exists a constant
$B_{\mathbb{K},m}^{\mathrm{mult}}\geq1$ such that for all $m$--linear forms
$T:\ell_{\infty}^{n}\times\cdots\times\ell_{\infty}^{n}\rightarrow\mathbb{K}$
and all positive integers $n$,
\begin{equation}
\left(  \sum_{j_{1},...,j_{m}=1}^{n}\left\vert T(e_{j_{1}},...,e_{j_{m}%
})\right\vert ^{\frac{2m}{m+1}}\right)  ^{\frac{m+1}{2m}}\leq B_{\mathbb{K}%
,m}^{\mathrm{mult}}\left\Vert T\right\Vert . \label{u88}%
\end{equation}

\bigskip These inequalities are, in some sense, predecessors of the
multilinear theory of absolutely summing operators (for details and recent
results on multilinear summing operators we refer to \cite{botelho, matos,
popa, pilar} and the references therein).

\bigskip The achievement of the optimal values of $B_{\mathbb{K}%
,m}^{\mathrm{mult}}$ and/or $C_{m,p}^{\mathbb{K}}$ is a quite challenging
problem and seems to be far from a definitive answer. From \cite{bohr} we know
that $B_{\mathbb{K},m}^{\mathrm{mult}}$ are dominated by constants with
sublinear growth, or more specifically,%
\begin{align*}
B_{\mathbb{C},m}^{\mathrm{mult}} &  <1.3\cdot m^{\frac{2-\log2-\gamma}{2}%
}<1.3\cdot m^{0.36482},\\
B_{\mathbb{R},m}^{\mathrm{mult}} &  <m^{\frac{1-\gamma}{2}}<m^{0.21139},
\end{align*}
where $\gamma$ from now on denotes the Euler--Mascheroni constant. On the
other hand, the best known upper bounds for the constants in (\ref{i99}) were,
until very recently, $\left(  \sqrt{2}\right)  ^{m-1}$ (see \cite{alb,
dimant}). In 2014, the upper estimate $\left(  \sqrt{2}\right)  ^{m-1}$ was
improved in the papers \cite{ARAUJO, DIOGO}:

\begin{itemize}
\item (Araujo \textit{et al}. \cite{DIOGO}) For all integers $m\geq2$ and
$p\in\left[  2m,\infty\right)  $ we have
\begin{align}
C_{m,p}^{\mathbb{R}}  &  <\left(  \sqrt{2}\right)  ^{\frac{2m\left(
m-1\right)  }{p}}\left(  1.3\cdot m^{0.36482}\right)  ^{\frac{p-2m}{p}%
},\label{yu9}\\
C_{m,p}^{\mathbb{C}}  &  <\left(  \frac{2}{\sqrt{\pi}}\right)  ^{\frac
{2m(m-1)}{p}}\left(  m^{0.21139}\right)  ^{\frac{p-2m}{p}}.\nonumber
\end{align}

\item (Araujo \textit{et al}. \cite{ARAUJO}) For all integers $m\geq2$ and
$p\in\left(  2m^{3}-4m^{2}+2m,\infty\right)  $ we have
\begin{align}
C_{m,p}^{\mathbb{R}}  &  <1.3\cdot m^{\frac{2-\log2-\gamma}{2}}<1.3\cdot
m^{0.36482}\label{yu10}\\
C_{m,p}^{\mathbb{C}}  &  <m^{\frac{1-\gamma}{2}}<m^{0.21139}.\nonumber
\end{align}

\end{itemize}

\bigskip A close look at (\ref{yu9}) and (\ref{yu10}) shows a surprising lack
of continuity when $p=2m^{3}-4m^{2}+2m$, i.e., by making $p\rightarrow
2m^{3}-4m^{2}+2m$ in (\ref{yu10}) we do not recover the estimate (\ref{yu9})
for $p=2m^{3}-4m^{2}+2m$. In this note we provide better estimates for the
case $2m\leq p\leq$ $2m^{3}-4m^{2}+2m$, improving (\ref{yu9}) whenever
$m\geq3$ (when $m=2$ our argument would provide the same constants from
\cite{ARAUJO}). Moreover, our estimates are continuous when compared with
(\ref{yu10}). More precisely, we prove the following result:

\begin{theorem}
\label{65l}\bigskip Let $m\geq3$ be a positive integer. For $2m\leq p\leq$
$2m^{3}-4m^{2}+2m$ there is a constant $C_{m,p}^{\mathbb{K}}\geq1$ such that,
for all $m$--linear forms $T:\ell_{p}^{n}\times\cdots\times\ell_{p}%
^{n}\rightarrow\mathbb{K}$ and all positive integers $n$,%
\[
\left(  \sum_{j_{1},...,j_{m}=1}^{n}\left\vert T(e_{j_{1}},...,e_{j_{m}%
})\right\vert ^{\frac{2mp}{mp+p-2m}}\right)  ^{\frac{mp+p-2m}{2mp}}\leq
C_{m,p}^{\mathbb{K}}\left\Vert T\right\Vert ,
\]
and
\end{theorem}

\begin{align}
C_{m,p}^{\mathbb{R}}  &  <\left(  1.3\cdot m^{\frac{2-\log2-\gamma}{2}%
}\right)  ^{\left(  m-1\right)  \left(  \frac{2m-p+mp-2m^{2}}{m^{2}%
p-2mp}\right)  }\left(  \sqrt{2}\right)  ^{\frac{1}{mp\left(  m-2\right)
}\left(  p-2m-mp+6m^{2}-6m^{3}+2m^{4}\right)  },\label{yhb}\\
C_{m,p}^{\mathbb{C}}  &  <\left(  m^{\frac{1-\gamma}{2}}\right)  ^{\left(
m-1\right)  \left(  \frac{2m-p+mp-2m^{2}}{m^{2}p-2mp}\right)  }\left(
\frac{2}{\sqrt{\pi}}\right)  ^{\frac{1}{mp\left(  m-2\right)  }\left(
p-2m-mp+6m^{2}-6m^{3}+2m^{4}\right)  },\nonumber
\end{align}
where $\gamma$ is the Euler--Mascheroni constant.

\begin{remark}
When $p=$ $2m^{3}-4m^{2}+2m$ we easily note that the estimates (\ref{yhb})
coincide with (\ref{yu10}). It is also not difficult to verify that the
estimates of (\ref{yhb}) are better than (\ref{yu9}).
\end{remark}

\section{The proof}

The proof follows the lines of the proof of \cite{DIOGO} with a technical
change in the interpolation procedure. Let
\[
s=\frac{2mp}{mp+p-2m}%
\]
and note that $\frac{2m}{m+1}\leq s\leq2.$ Let also%
\[
\lambda_{0}=\frac{2s}{ms+s-2m+2}%
\]
and note that
\[
\lambda_{0}<s\leq2.
\]
Since%
\[
\frac{m-1}{s}+\frac{1}{\lambda_{0}}=\frac{m+1}{2},
\]
from the generalized Bohnenblust--Hille inequality (see \cite[Theorem
1.1]{alb}) we know that there is a constant $C_{m}\geq1$ such that for all
$m$-linear forms $T:\ell_{\infty}^{n}\times\cdots\times\ell_{\infty}%
^{n}\rightarrow\mathbb{K}$ we have, for all $i=1,....,m$ and all positive
integers $n,$%

\begin{equation}
\left(  \sum\limits_{j_{i}=1}^{n}\left(  \sum\limits_{\widehat{j_{i}}=1}%
^{n}\left\vert T\left(  e_{j_{1}},...,e_{j_{m}}\right)  \right\vert
^{s}\right)  ^{\frac{1}{s}\lambda_{0}}\right)  ^{\frac{1}{\lambda_{0}}}\leq
C_{m}\left\Vert T\right\Vert . \label{78}%
\end{equation}
Above, as usual, $\sum\limits_{\widehat{j_{i}}=1}^{n}$ means that we are
summing over all $j_{k}$ for all $k\neq i.$

Now we need a good estimate for the constant $C_{m}$ in (\ref{78}). For this
task, note that the multiple exponent%
\[
\left(  \lambda_{0},s,s,...,s\right)
\]
is the result of the interpolation (in the sense of \cite[Section 2]{alb}) of
the multiple exponents%

\[
\left\{
\begin{array}
[c]{c}%
E_{1}=\left(  \frac{2m-2}{m},....,\frac{2m-2}{m},2\right)  \\
E_{2}=\left(  \frac{2m-2}{m},\frac{2m-2}{m},....,2,\frac{2m-2}{m}\right)  \\
\vdots\\
E_{m-1}=\left(  \frac{2m-2}{m},2,\frac{2m-2}{m},....,\frac{2m-2}{m}\right)  \\
E_{m}=\left(  1,2,....,2\right)
\end{array}
\right.
\]
with, respectively,
\[%
\begin{array}
[c]{c}%
\theta_{1}=\cdots=\theta_{m-1}=\frac{2m-p+mp-2m^{2}}{m^{2}p-2mp}\\
\theta_{m}=1-\left(  m-1\right)  \left(  \frac{2m-p+mp-2m^{2}}{m^{2}%
p-2mp}\right)  .
\end{array}
\]
Above note that we have used the hypothesis $2m\leq p\leq$ $2m^{3}-4m^{2}+2m.$
From \cite{ARAUJO, bohr} we know that for the multiple exponents
$E_{1},...,E_{m-1}$ the inequality (\ref{78}) is valid with constants
$1.3\cdot m^{\frac{2-\log2-\gamma}{2}}$ for real scalars and $m^{\frac
{1-\gamma}{2}}$ for complex scalars. On the other hand, it is well known that
the constant associated with the multiple exponent $\left(  1,2,...,2\right)
$ is $\left(  \sqrt{2}\right)  ^{m-1}$ for real scalars (see \cite{jnt}) and
less than or equal to $\left(  \frac{2}{\sqrt{\pi}}\right)  ^{m-1}$ for
complex scalars$.$ Therefore, the optimal constant associated to the multiple
exponent
\[
\left(  \lambda_{0},s,s,...,s\right)
\]
is less than or equal to \
\[
\left(  1.3\cdot m^{\frac{2-\log2-\gamma}{2}}\right)  ^{\left(  m-1\right)
\left(  \frac{2m-p+mp-2m^{2}}{m^{2}p-2mp}\right)  }\left(  \sqrt{2}\right)
^{\frac{1}{mp\left(  m-2\right)  }\left(  p-2m-mp+6m^{2}-6m^{3}+2m^{4}\right)
}%
\]
for real scalars and%
\[
\left(  m^{\frac{1-\gamma}{2}}\right)  ^{\left(  m-1\right)  \left(
\frac{2m-p+mp-2m^{2}}{m^{2}p-2mp}\right)  }\left(  \frac{2}{\sqrt{\pi}%
}\right)  ^{\frac{1}{mp\left(  m-2\right)  }\left(  p-2m-mp+6m^{2}%
-6m^{3}+2m^{4}\right)  }%
\]
for complex scalars.

Let%
\[
\lambda_{j}=\frac{\lambda_{0}p}{p-\lambda_{0}j}%
\]
for all $j=1,....,m.$ Note that%
\[
\lambda_{m}=s
\]
and that the conjugate number of $\left(  \frac{p}{\lambda_{j}}\right)  $ is
$\frac{\lambda_{j+1}}{\lambda_{j}}$, i.e.,
\[
\left(  \frac{p}{\lambda_{j}}\right)  ^{\ast}=\frac{\lambda_{j+1}}{\lambda
_{j}}%
\]
for all $j=0,...,m-1$. Now the proof follows straightforwardly the lines of
\cite{DIOGO}.


\begin{figure}[h]
\centering
\includegraphics[scale=0.6]{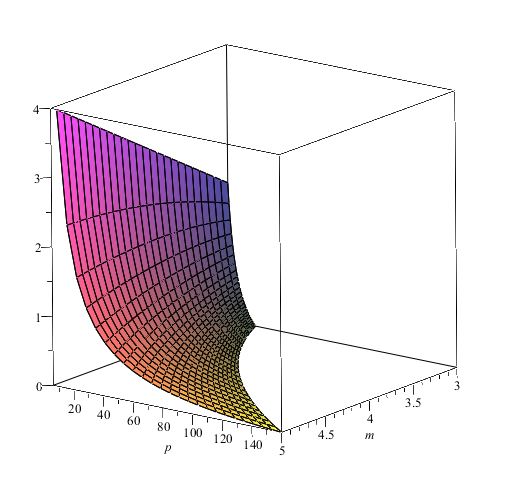}\caption{Exponent of $\sqrt{2}$ in
(\ref{yhb})}%
\label{figure2}%
\end{figure}

\section{Final remark: a more general approach}

The Hardy--Littlewood inequality for multilinear forms has the following
version for multiple exponents:

\textbf{Theorem (Generalized Hardy--Littlewood inequality, (\cite{alb},
2014)). }\textit{Let }$m\geq2$ be a positive integer, $2m\leq p\leq\infty$
\textit{and} $\mathbf{q}:=(q_{1},...,q_{m})\in\left[  \frac{p}{p-m},2\right]
^{m}$ be \textit{such that}
\begin{equation}
\frac{1}{q_{1}}+\cdots+\frac{1}{q_{m}}=\frac{mp+p-2m}{2p}. \label{ppoo}%
\end{equation}
\textit{Then there exists a constant} $C_{m,p,\mathbf{q}}^{\mathbb{K}}\geq1$
\textit{such that}
\begin{equation}
\left(  \sum_{j_{1}=1}^{n}\left(  \sum_{j_{2}=1}^{n}\left(  \cdots\left(
\sum_{j_{m}=1}^{n}\left\vert T(e_{j_{1}},...,e_{j_{m}})\right\vert ^{q_{m}%
}\right)  ^{\frac{q_{m-1}}{q_{m}}}\cdots\right)  ^{\frac{q_{2}}{q_{3}}%
}\right)  ^{\frac{q_{1}}{q_{2}}}\right)  ^{\frac{1}{q_{1}}}\leq
C_{m,p,\mathbf{q}}^{\mathbb{K}}\left\Vert T\right\Vert \label{g33}%
\end{equation}
\textit{for all }$m$\textit{--linear forms }$T:\ell_{p}^{n}\times\cdots
\times\ell_{p}^{n}\rightarrow\mathbb{K}$\textit{\ and all positive integers
}$n.$\textit{\ }

\medskip In \cite{ARAUJO} the following result is proved:

\begin{theorem}
\label{999}(\cite{ARAUJO})Let $m\geq2$ be a positive integer and
$2m<p\leq\infty$. Let also $\mathbf{q}:=\left(  q_{1},...,q_{m}\right)
\in\left[  \frac{p}{p-m},2\right]  ^{m}$ be such that $\frac{1}{q_{1}}%
+\cdots+\frac{1}{q_{m}}=\frac{mp+p-2m}{2p}$. If $\max q_{i}<\frac{2m^{2}%
-4m+2}{m^{2}-m-1}$, then%
\begin{align*}
C_{m,p,\mathbf{q}}^{\mathbb{R}} &  <1.3\cdot m^{\frac{2-\log2-\gamma}{2}%
}<1.3\cdot m^{0.36482},\\
C_{m,p,\mathbf{q}}^{\mathbb{C}} &  <m^{\frac{1-\gamma}{2}}<m^{0.21139}.
\end{align*}
\bigskip
\end{theorem}

Using the same ideas of the previous section and following the lines of
\cite{ARAUJO} we can prove the following result (the proof is long, although
the arguments are similar to the previous, and we left the details for the
interested reader):

\begin{theorem}
\label{765}Let $m\geq2$ be a positive integer and $2m\leq p\leq\infty$. Let
also $\mathbf{q}:=\left(  q_{1},...,q_{m}\right)  \in\left[  \frac{p}%
{p-m},2\right]  ^{m}$ be such that
\[
\frac{1}{q_{1}}+\cdots+\frac{1}{q_{m}}=\frac{mp+p-2m}{2p}%
\]
and
\begin{equation}
\max q_{i}\geq\frac{2m^{2}-4m+2}{m^{2}-m-1}. \label{7776}%
\end{equation}
\bigskip For all $m$--linear forms $T:\ell_{p}^{n}\times\cdots\times\ell
_{p}^{n}\rightarrow\mathbb{K}$ and all positive integers $n$, we have%

\[
\left(  \sum_{j_{1}=1}^{n}\left(  \sum_{j_{2}=1}^{n}\left(  \cdots\left(
\sum_{j_{m}=1}^{n}\left\vert T(e_{j_{1}},...,e_{j_{m}})\right\vert ^{q_{m}%
}\right)  ^{\frac{q_{m-1}}{q_{m}}}\cdots\right)  ^{\frac{q_{2}}{q_{3}}%
}\right)  ^{\frac{q_{1}}{q_{2}}}\right)  ^{\frac{1}{q_{1}}}\leq
C_{m,p,\mathbf{q}}^{\mathbb{K}}\left\Vert T\right\Vert ,
\]
with%
\begin{align*}
C_{m,p,\mathbf{q}}^{\mathbb{C}}  &  \leq\left(  \frac{2}{\sqrt{\pi}}\right)
^{\left(  m-1\right)  \theta_{1}}\left(  \eta_{\mathbb{C},m}\right)
^{\theta_{2}},\\
C_{m,p,\mathbf{q}}^{\mathbb{R}}  &  \leq\left(  \sqrt{2}\right)  ^{\left(
m-1\right)  \theta_{1}}\left(  \eta_{\mathbb{R},m}\right)  ^{\theta_{2}},
\end{align*}

where $\left(  \theta_{1},\theta_{2}\right)  =\left(  1,0\right)  $ if $m=2$
and
\begin{equation}
\theta_{1}=1-\frac{\left(  m+1\right)  \left(  2-\max q_{i}\right)  \left(
m-1\right)  ^{2}}{\left(  m^{2}-m-2\right)  \max q_{i}}\text{ and }\theta
_{2}=\frac{\left(  m+1\right)  \left(  2-\max q_{i}\right)  \left(
m-1\right)  ^{2}}{\left(  m^{2}-m-2\right)  \max q_{i}} \label{teta2}%
\end{equation}
for $m\geq3.$
\end{theorem}

\bigskip

Note that now we have continuity between the estimates of Theorem \ref{765}
and Theorem \ref{999} since
\[
\frac{\left(  m+1\right)  \left(  2-\max q_{i}\right)  \left(  m-1\right)
^{2}}{\left(  m^{2}-m-2\right)  \max q_{i}}=1
\]
when $\max q_{i}=\frac{2m^{2}-4m+2}{m^{2}-m-1}.$

\end{document}